% Critical metrics for the zeta function of the hodge Laplacian
% Kate Okikiolu and Caitlin Wang

\input amstex.tex

\documentstyle{amsppt}
\magnification \magstephalf
\pagewidth{6.2in}
\pageheight{8.0in}
%\hcorrection{0.6in}
\NoRunningHeads
\def\makefootline{\baselineskip=49pt \line{\the\footline}}
\define\hx{h\!\cdot\!\xi}
\define\xhx{\xi\!\cdot\!h\!\cdot\!\xi}
\define\s{\sigma}
\define\la{\lambda}
\define\al{\alpha}

\define\tr{\operatorname{tr}}

\define\trace{\operatorname{trace}}
\define\Hess{\operatorname{Hess}}
\define\End{\operatorname{End}}

\define\sgn{\operatorname{sgn}}
\define\llangle{{\langle\!\langle}}
\define\rrangle{{\rangle\!\rangle}}
\define\lm{ \left( \matrix }
\define\mr{\endmatrix \right)}
\def\today{\ifcase\month\or
    January\or February\or March\or April\or May\or June\or
    July\or August\or September\or October\or November\or December\fi
    \space\number\day, \number\year}

\topmatter
\title
HESSIAN OF THE ZETA FUNCTION FOR THE LAPLACIAN ON FORMS
\endtitle
\author  K. Okikiolu and C. Wang
\endauthor
\thanks  The first author was supported  by the National
Science Foundation  \#DMS-9703329.  \endthanks \abstract  Let $M$
be a compact closed $n$-dimensional manifold. Given a Riemannian
metric on $M$,  we consider the zeta functions $Z(s)$ for the de
Rham Laplacian and the Bochner Laplacian on $p$-forms.  The
hessian of $Z(s)$ with respect to variations of the metric is
given by a pseudodifferential operator $T_s$. When the real part
of $s$ is less than $n/2-1$ we compute the principal symbol of
$T_s$. This can be used to determine whether a general critical
metric for $(d/ds)^k Z(s)$ has finite index, or whether it is an
essential saddle point.
\endabstract
\endtopmatter

\noindent{\bf 1. Introduction.}
\medskip

Over the last fifteen years, several results have been proved
which identify  maximal or minimal metrics  for the determinant of
a geometric elliptic operator on a compact manifold. In [OPS1],
surfaces of constant curvature were shown to maximize the
determinant of the Laplacian within the conformal class.
 Most  higher dimensional cases studied
 involve extremizing a conformally covariant operator
 within a conformal class.
   In these cases local variation formulas are available, [B\O], [Po], and
 global results have been obtained, see  [BCY], [Br],
[CY].  (Some results also exist for manifolds which are
 non-closed or non-compact, see for example [BFKM], [CQ], [HZ], [OPS2], [OPS3].)
In addition, extremal metrics for special values of the zeta
function for conformally covariant operators have been studied,
see [Mo].  The zeta function for the de Rham Laplacian on forms is
interesting because of its connection to topology.
 For the Laplace-Beltrami operator $\Delta_0$,
the standard $3$-sphere is a local maximum for $\det\Delta_0$, but
if $n\geq 7$ is odd, the standard $n$-sphere is a saddle for
$\det\Delta_0$, see [Ok1], [Ri].  Beyond the scalar case, the
determinant of the de Rham Laplacian was studied in relation to
analytic torsion, [Ch], [M\"{u}], [RS] and  the variation of the
zeta function was studied in [Ro], but there does not seem to have
been progress on the subject of critical metrics since then. (In
the non-scalar case there is not a great deal known about the
spectrum of the Laplacian on forms. See [Ba], [BGV],  [CT], [CC],
[Fe], [Do], [Go], [Lo], [Ta] for results and further references.)
In [Ok2] it was shown that for general geometric operators of
Laplace type, the Hessian of the zeta function is given by a
pseudodifferential operator $T_s$, and when $\Re s<n/2-1$, a
formula for the principal symbol $u_s(x,\xi)$ was given. In the
current paper we compute $u_s$ explicitly for the case of the de
Rham and Bochner Laplacians on $p$-forms. This can be applied to
compute the behavior of the determinant and zeta function in the
neighborhood of critical metrics. The computation of $u_s$ is
somewhat lengthy for the de Rham Laplacian. The proof given here
can serve as a blueprint for carrying out the calculation for any
operator of Laplace type. Happily, the results of the calculations
are simple to state, see Theorem 1.

In order to present the results we set up some notation. Let $M$
be a closed, compact $n$-dimensional manifold, let $\chi$ be a
character of the fundamental group of $M$,  let $L_\chi$ denote
the complex line over $M$  twisted by $\chi$, let
$\Omega^p(M,L_\chi)$ denote the space of smooth $p$-forms valued
in $L_\chi$, and let
 $d$ denote the  exterior differential.
  For the Riemannian metric $g$ on $M$, let $dV$ be the volume form on $M$
  and let $V$ be the volume of $M$. For any tensor bundle $E$ over $M$
  with values in $L_\chi$, the metric $g$ naturally gives rise
  to  an inner product $\langle\ ,\ \rangle_g$ on $E$, an an inner product
  $\llangle\ ,\ \rrangle_g$ on the smooth sections of $E$, and a
   connection $\nabla$ on $E$.
 We consider the de Rham Laplacian $d^* d+d d^*$ and the Bochner
  Laplacian $\nabla^*\nabla$ on $\Omega^p(M,L_\chi)$.  These differ
  by a curvature operator.  We use the notation $\Delta_p$
  to denote either of the two operators.
Let
$$
\la_1\ \leq\ \la_2\ \leq\ \dots
$$
be the non-zero eigenvalues of $\Delta_p$. The spectral zeta
function
$$
Z(s)\ =\ Z(\Delta_p)(s)\ =\ \sum_{j=1}^\infty \la_j^{-s}
$$
converges for $\Re s>n/2$ and extends to a meromorphic function
with simple poles located at
$$
s\ =\ \cases n/2,\ n/2-1,\ \dots, \ 1,  \qquad & n\
\text{even},\\
n/2,\ n/2-1,\ \dots\ \dots,\ \qquad & n\ \text{odd}.\endcases
$$
The spectral determinant of $\Delta_p$ is defined to be
$$
\det \Delta_p\ =\ e^{-Z^\prime(0)}.
$$
Define the modified zeta function to be
$$
\Cal Z(s)\ =\  \Cal Z(\Delta_p)(s)\ =\
\frac{\Gamma(s)Z(s)}{\Gamma(s-n/2)}
 \ +\ \frac{\dim\ker \Delta_p}{s\ \Gamma(s-n/2)}    ,\tag 1.1
 $$
which is entire in $s$. We remark that the dimension,
$\beta_p(\chi)$,  of the kernel of
 the de Rham Laplacian is a topological invariant, and is
 consequently invariant under changes of the metric.

 Let $S^2M$ denote the bundle of symmetric tensor fields of type $(0,2)$
 on $M$, and let $C^\infty(S^2M)$ denote its space of smooth
 sections.
 The Hessian of the modified zeta function is the bilinear form on
$C^\infty(S^2M)$ defined by
$$
\Hess \Cal Z(s)(h,h)\ =\ \frac{d^2}{d\al^2}\biggl|_{\al=0} \Cal
Z(\Delta_p(g+\al h))(s).
$$

We recall the following result which holds for the general
geometric operator $\Delta$ of Laplace type defined on a tensor
bundle.
 \proclaim{[Ok2, Theorem 1]}  For $s\in \Bbb C$,  there exists a unique
symmetric pseudodifferential operator $T_s=T_s(\Delta
):C^\infty(S^2M)\to C^\infty(S^2M)$ such that
$$
\Hess \Cal Z(s)(h,h)\ =\ \llangle h, T_s h\rrangle_g. \tag 1.2
$$
The operator $T_s$ is analytic in $s$.  For  $s\notin n/2+\Bbb
N^+$, there exist polyhomogeneous pseudodifferential operators
$U_s$ and $V_s$ of degrees $n-2s$ and $2$ respectively such that
$T_s=U_s+V_s$.  The operators $U_s$ and $V_s$ are meromorphic in
$s$ with simple poles located in $ n/2+\Bbb N^+$. (The poles of
$U_s$ and $V_s$ cancel  in the sum $T_s=U_s+V_s$, but the symbol
expansion of $T_s$ for $s\in n/2+\Bbb N^+$ involves logarithmic
terms.) For general $s$, the polyhomogeneous symbol expansion of
$U_s$ is computable from the complete symbol of the operator
$\Delta$. In particular, there is a simple algorithm to compute
the term $u_s$ of homogeneity $n-2s$.
 Furthermore, we can differentiate in
$s$ to obtain
$$
\Hess (d/ds)^k \Cal Z(s)(h,h)\ =\ \llangle h, (d/ds)^k T_s
h\rrangle_g,\tag 1.3
$$
and the principal symbol of $(d/ds)^k U_s$ is equal to the leading
order  term of  $(d/ds)^k u_s$  (provided this does not vanish
identically).
\endproclaim

The purpose of this paper is to compute the symbol $u_s(x,\xi)$
for the de Rham and Bochner Laplacians.  This is the principal
symbol of $T_s$ when $\Re s<n/2-1$, (providing it does not vanish
identically) and so for $s$ in this range, $u_s$ governs to some
extent how $(d/ds)^k Z(s)$ behaves in the neighborhood of a
critical metric. In particular, $u_s(x,\xi)$ generally governs the
behavior of $\log\det\Delta$.  We remark that  the above theorem
was stated in [Ok2] for the case when the Laplacian  has been
normalized to be scale invariant, but this assumption is easily
seen to be unnecessary.

For a point $x\in M$, we write $\End(T_x^*M)$ for the
endomorphisms of the cotangent space $T_x^*M$, and we write
``$\tr$" for the  trace function on $\End(T_x^*M)$. It is
convenient to introduce the variable $S=s-n/2$, and set
 $$
 C= C(n,s)\ =\ \left(\frac1{4\pi}\right)^{n/2}
\frac{\Gamma(-S+1)^2}{\Gamma(-2S+2)}. \tag 1.4
$$
\proclaim{Theorem 1} Suppose $g$ is a metric on $M$, $x$ is a
point of $M$, and $h\in S^2_xM$.  Set
$$
H\ =\ hg^{-1} \ \in\ \End(T_x^*M).\tag 1.5
$$
For a non-zero  covector $\xi\in T_x^*M$, let
$$
\Pi_\xi: T^*_x(M)\ \to\ T^*_x(M)
$$
denote the $g$-orthogonal projection onto the space spanned by
$\xi$, and let
$$
\Pi_\xi^\perp\ =\ I\ -\ \Pi_\xi. \tag 1.6
$$
Let $\Delta_p$ be the de Rham or Bochner Laplacian on $p$-forms,
and in the case of the Bochner Laplacian assume that the dimension
of its kernel has constant dimension when $g$ is deformed.  Then
at the metric $g$,
$$
  \langle h,u_s(x,\xi) h\rangle_g\ =\
  2C(n,s)|\xi|^{n-2s}\left(\ f_1(\Delta_p)\ \tr\ \left(( H\Pi_\xi^\perp)^2\right)
\ +\
  f_2(\Delta_p)\ \left(\,\tr\, H\,\Pi_\xi^\perp\,\right)^2 \right),
  \tag 1.7
$$
where
$$
\align
 &f_1(\nabla^*\nabla)(n,p,s)\ =\ 4\lm n-2\\ p-1\mr (S-1/2)\ +\ \frac12\lm n\\
p\mr ,\tag 1.8\\
 &f_2(\nabla^*\nabla)(n,p,s)\ =\ \lm n\\ p\mr(S^2+S-3/4)\ +\ \frac14\lm
 n\\ p\mr .\tag 1.9\\
&f_1(d^*d+dd^*)(n,p,s)\ =\ 4\lm n-2\\p-1\mr (S^2+S-3/4)\ +\
\frac12\lm
n\\p\mr, \tag 1.10\\
 &f_2(d^*d+dd^*)(n,p,s)\ =\ \left( \lm n\\p \mr -4\lm n-2\\p-1\mr\right)
 (S^2+S-3/4)
 \ +\ \frac14\lm n\\p\mr , \tag 1.11\\
\endalign
$$
\endproclaim
\medskip

Before discussing applications, we point out some features of the
result.  Firstly one notices that for  $\xi,\eta\in T_x^*M$, if
$h\in S^2_xM$ is given by
 $$
 h\ =\ \xi\otimes \eta\ +\ \eta\otimes \xi
 \tag 1.12
 $$
then
$$
\Pi_\xi^\perp\ H\ \Pi_\xi^\perp\ =\ 0.
$$
Hence for any $k\in S^2_xM$, writing $K=kg^{-1}$ we have
$$
  \langle k,u_s(x,\xi) h\rangle_g\ =\
  2C(n,s)|\xi|^{n-2s} \left(\ f_1\trace\ (K\Pi_\xi^\perp H\Pi_\xi^\perp)
  \ +\ f_2  \left(\,\trace\, K\,\Pi_\xi^\perp\,\right)
  \left(\,\trace\,  H\,\Pi_\xi^\perp\,\right)\ \right)
  \ \  =\ \ 0. \tag 1.13
$$
This reflects on a micro-local level the fact that if $\phi(\al)$
is a one parameter family of diffeomorphisms of $M$ with $\phi(0)$
equal to the identity and we set
$$
h\ =\ \frac{d\phi(\al)^* g}{d\al}\biggl|_{\al=0},
$$
then  $\Hess \Cal Z(s)(k,h)$ vanishes for every $k\in
C^\infty(S^2M)$.  It is easy to see that if $h\in S^2M$ does not
have the form in (1.12) then
$$
\trace\ \left(( H\Pi_\xi^\perp)^2\right)\ >\ 0,\qquad\qquad\qquad
\frac1{n-1}\trace\left(( H\Pi_\xi^\perp)^2\right)\ \geq\
\left(\,\trace\, H\,\Pi_\xi^\perp\,\right)^2\ \geq\ 0. \tag 1.14
$$
Secondly, we remark that instead of considering the full de Rham
Laplacian, we can consider the
 operator $d^* d$ acting on $p$-forms.  The zeta function
 for this operator is
 $$
 Z((d^* d)_p)\, (s)\ =\ \sum_{q=0}^p (-1)^{p-q}\ Z((d^* d+dd^*)_q)(s).
 $$
 From Theorem 1, we find that the symbol $u_s(x,\xi)$ for
 the Hessian of $\Cal Z_{d^* d} (s)$ is given by
 $$
 \langle h,u_s(x,\xi) h\rangle_g
 \ =\ C(n,s)\left(\ f_1(n-1,p,s)\tr\ \left(( H\Pi_\xi^\perp)^2\right)\ +\
 f_2(n-1,p,s)\left(\,\tr\, H\,\Pi_\xi^\perp\,\right)^2
 \ \right)|\xi|^{n-2s}. \tag 1.15
 $$
Thirdly we remark that the analytic torsion $\tau$  defined in the
acyclic case by
$$
\frac12\log\tau\ =\ \sum_{p=0}^n (-1)^{p+1} p\ \log\det
(d^*d+dd^*)_p
$$
is a topological invariant.  This implies that the corresponding
alternating sum of symbols $u_0(x,\xi)$ for $d^*d+dd^*$ on $p$
forms must vanish, and this is easily seen to be the case. In fact
a stronger result is easily seen to be true from Theorem 1, namely
for both the de Rham and the Bochner Laplacian $\Delta_p$ and for
all $s$ and all  $k=0,1,\dots$,
$$
\sum_{p=0}^n (-1)^{p+1} p^k\ {u_s}(\Delta_p)(x,\xi)\ =\ 0. \tag
1.16
$$
Lastly we remark that the character $\chi$ plays no role in our
results which come from local calculations.  The only thing one
sees in the results related to the line bundle $L_\chi$ is the
factor ``$2$" in the right hand side of (1.7) which comes from the
fact that the complex line has real dimension $2$.

In a subsequent paper, we will show how Theorem 1 implies
following results.

\proclaim{Theorem 2}  Suppose $n\geq 3$, $s$ is real and
$s<n/2-1$, and $\Delta_p$ is either the de Rham or the Bochner
Laplacian. For $j=1,2$,  we consider the real numbers
$f_j=f_j(\Delta_p)(n,p,s)$.

\noindent (a). If $f_1$ and $f_1/(n-1)+f_2$ and are both positive
or both negative, then  every critical metric $g$ for
 $$
 (-1)^k \sgn(f_1)\ (d/ds)^k \Cal Z(s)
 \tag 1.17
 $$
  has finite index.  (That is, there exists a submanifold
$\Cal X$ of $\Cal M$ of finite codimension, containing $g$, such
$(-1)^k \sgn(f_1)\, (d/ds)^k \Cal Z(s)$ has a local minimum at $g$
on $\Cal X$.)

\noindent (b). Conversely, if $f_1$ and $f_1/(n-1)+f_2$ have
opposite signs, then every critical metric $g$ for $(d/ds)^k \Cal
Z(s)$ is an essential saddle point. (That is, there are two
infinite dimensional submanifolds $\Cal X,\Cal Y$, of $\Cal M$,
both containing $g$, such that $(d/ds)^k\Cal Z(s)$ has a local
minimum at $g$ on $\Cal X$ and a local maximum at $g$ on $\Cal
Y$.)
\endproclaim

\proclaim{Corollary} Let $n\geq 3$. If $n=1,\mod 4$, every
critical metric for $\det(d^*d+dd^*)_p$ has finite index.  If
$n=3,\mod 4$, every critical metric for $1/\det(d^*d+dd^*)_p$ has
finite index. The same is true if the de Rham Laplacian is
replaced by the operator $d^* d$ or the operator $dd^*$ acting on
$p$ forms. However, for the Bochner Laplacian on forms of degree
$p\neq 0,n$, the determinant may have an essential saddle point.
The lowest dimension in which this  happens is $n=4$.
\endproclaim

We remark that because we have not normalized the Laplacian to be
scale invariant, a critical metric in the results above is one
which is critical for the given functional under variations of the
metric which preserve the total volume.

 The rest of the paper is devoted to proving
Theorem 1. We work throughout in the real case,  where the
Laplacian is acting on real valued $p$-forms.  Because the
calculation is entirely local, the character $\chi$ plays no role.
The effect of going from the real case to the case to the case of
the complex line $L_\chi$ is just to multiply $u_s$ by the factor
$2$, and this is the factor $2$ which appears in (1.7). In Section
2, we discuss the general formula for $u_s$ which was given in
[Ok2], and we rewrite it in terms of the {\it coefficient symbols}
of the first variation of the Laplacian. In Section 3, we compute
the leading order terms in the first variation of $\Delta_p$,
where $\Delta_p$ is the Bochner or de Rham Laplacian. This is
carried out in coordinates utilizing (3.7) as suggested by a
referee for an early version of the paper which dealt only with
the de Rham operator. Our change in notation has the advantage
that the Bochner and de Rham Laplacians can be dealt with in the
same way, and it will also hopefully make the proofs easier to
follow.   In Section 4, we identify the matrix entries of the
first variation of the Laplacian.  In Sections 5 and 6, we
substitute these matrix entries into our formula for $u_s$, thus
proving Theorem 1.  For case of the de Rham Laplacian, it is
essential to work in coordinates  $\xi$ in which $h$ is diagonal,
since otherwise the calculation becomes extremely lengthy. It is
justified to choose to work in such coordinates because $u_s$ is
the principal symbol of a pseudodifferential operator, and hence
does not depend on the choice of coordinates. Section 7 contains
some observations on the calculation of $u_s$ for the de Rham
Laplacian in general coordinates.

\medskip
\medskip
%%%%%%%%%%%%%%%%%%%%%%%%%%%%%%%%%%%%%%%%%%%%%%%%%%%%%%%%%%%%%%%%%%%
\noindent{\bf 2. The general formula for $u_s$.} \medskip

Our starting point to prove Theorem 1, is the following result
which gives $u_s$ in the general case.

 \proclaim{[Ok2, Theorem 2]} The symbol
$u_s$ for the general geometric operator $\Delta$ of Laplace type
on the tensor bundle $E$ of real dimension $N$ can be computed as
follows. Write $\Delta_p '=(d/d\s)|_{\s=0} \Delta_p (g+\s h)$ and
 let $x\in M$.  Take coordinates on $M$ which are orthonormal in the metric $g$ at
 the point $x$, and take a local trivialization of $E$ on a neighborhood of $x$ to
 obtain coordinates for $E$.
  Suppose that in these coordinates,  the
operator $\Delta_p ^\prime$ is given at the point $x$ by
$$
\Delta_p ^\prime\ =\ \sum_{\al,\beta,i,j} A_{\al,\beta}^{ij} \
(\partial^\al_w h_{ij}(x))\ \partial_w^\beta, \tag  2.1
$$
where  $\al$ and $\beta$ are multiindices and $A_{\al,\beta}^{ij}$
is an $N\times N$ real valued matrix. Set
$$
C(n,s)\ =\ \left(\frac1{4\pi}\right)^{n/2}
\frac{\Gamma(-S+1)^2}{\Gamma(-2S+2)}. \tag 2.2
$$

Then at $(x,\xi)\in T^*M$, the value of $u_s(x,\xi)\in\End
(S^2M)_x$ is given by
$$
(u_s(x,\xi)h)_{ij}\ =\ C(n,s)\sum\Sb |\al|+|\beta|=2\\
|\gamma|+|\delta|=2 \\ k,\ell\endSb
\,u_s(\partial^\al,\partial^\beta,\partial^\gamma,\partial^\delta,x,\xi)
\ \trace(A_{\al,\beta}^{ij}(x)\,A_{\gamma,\delta}^{k\ell}(x))\
h_{k\ell}, \tag 2.3
$$
where writing $I$ for the identity operator on $E_x$, the terms
$u_s(\partial^\al,\partial^\beta,\partial^\gamma,\partial^\delta,x,\xi)$
are given as follows:
$$
u_s(\partial_j\partial_k, \, I,\, \partial_p\partial_q,\,I, \,
x,\xi) \ =\ 4(4S^2-1) \xi_j\xi_k\xi_p\xi_q\,|\xi|^{n-2s-4} .\tag
2.4
$$
$$
u_s(\partial_j, \, \partial_k,\, \partial_p\partial_q,\,I,\,
x,\xi) \ =\   -2(4S^2-1)\xi_j\xi_k\xi_p\xi_q |\xi|^{n-2s-4}.\tag
2.5
$$
$$
u_s(\partial_j, \, \partial_k,\, \partial_p,\, \partial_q,\,
x,\xi) \ =\   (4S^2+2S-2)\xi_j\xi_k\xi_p\xi_q |\xi|^{n-2s-4}\ -\
 (2S-1)\delta_{kq}\xi_j\xi_p |\xi|^{n-2s-2}.\tag 2.6
$$
$$
u_s(I, \, \partial_j\partial_k,\, \partial_p\partial_q,I,\, x,\xi)
\ =\  (4S^2-2S)\xi_j\xi_k\xi_p\xi_q |\xi|^{n-2s-4}\ +\
(2S-1)\delta_{jk}\xi_p\xi_q |\xi|^{n-2s-2}. \tag 2.7
$$
$$
\multline u_s(I, \, \partial_j\partial_k,\, \partial_p,\,
\partial_q,\, x,\xi) \  =\ -(2S^2+S-1)\xi_j\xi_k\xi_p\xi_q
|\xi|^{n-2s-4}
 \\ +\   (S-1/2)(-\delta_{jk}\xi_p\xi_q+
 \delta_{jq}\xi_k\xi_p+\delta_{kq}\xi_j\xi_p) |\xi|^{n-2s-2}
\endmultline
\tag 2.8
$$
$$
\align &u_s (I,  \,  \partial_j\partial_k, \,I,\,
\partial_p
\partial_q,\, x,\xi)
\  =  \tag 2.9 \\
 (S^2+ & S) \xi_j\xi_k\xi_p\xi_q|\xi|^{n-2s-4}\ +\
\frac12(S-1)\delta_{jk}\xi_p\xi_q|\xi|^{n-2s-2}
 \ +\ \frac12(S-1)\xi_j\xi_k\delta_{pq}|\xi|^{n-2s-2}\\
 -\   \frac{S}{2}& (\delta_{jp}\xi_k\xi_q+\delta_{kq}\xi_j\xi_p
 +\delta_{jq}\xi_k\xi_p+\delta_{kp}\xi_j\xi_q)|\xi|^{n-2s-2}\ +\
  \frac{1}{4}(\delta_{jk}\delta_{pq}+\delta_{jp}\delta_{kq}
  +\delta_{jq}\delta_{kp})|\xi|^{n-2s}.
\endalign
$$
\endproclaim
\medskip
It is convenient to introduce the following notation.
$$
 \hx=\sum_i h_{ij}\xi_j,\quad\qquad \xhx=\sum_{i,j}
h_{ij}\xi_i\xi_j,\quad\qquad
|h|=\left(\sum_{i,j}h_{ij}^2\right)^{1/2},\quad\qquad \tr h=\sum_i
h_{ii}.\tag 2.10
$$
Note that the matrix for $H\Pi_\xi^\perp$ is
$$
h_{ij}\ -\ |\xi|^{-2}\sum_k  h_{ik}\xi_k\xi_j,
$$
so
$$
\align & \tr \left( \left( H\Pi_\xi^\perp\right)^2\right)= |h|^2\
-\
2|\xi|^{-2}|\hx|^2\ +\ |\xi|^{-4}(\xhx)^2, \tag 2.11\\
 & \left(
\tr H\Pi_\xi^\perp\right)^2= |\xi|^{-4}(\xhx)^2\ -\
2|\xi|^{-2}(\xhx)( \tr h)\ +\ (\tr h)^2.
\endalign
$$

 We will now write $u_s$ as a sum
$$
u_s\ =\ C(n,s)\left(u_s^{(1)}\ +\ u_s^{(2)}\ +\ u_s^{(3)}\ +\
u_s^{(4)}\right). \tag 2.12
$$
To define the symbols $u_s^{(j)}$, choose coordinates which are
orthonormal for the metric $g$ at the point $x$.  We notice that
(2.3) only involves those terms in (2.1) which have {\it total
degree} $2$, where the  {\it total degree} of a term
$$
(\partial^\al_w h_{ij}(x))\ \partial_w^\beta \tag 2.13
$$
is $|\al|+|\beta|$.  We  write
$$
\Delta'\ \sim\ \sum_{i,j} h_{ij}\partial_{ij}^2\ +\
\sum_{i,j,k,\ell} \left( A_{k\ell}^{ij}\left(\partial_k
h_{ij}\right)\partial_\ell\ +\
B_{k\ell}^{ij}\left(\partial_{k\ell}^2 h_{ij}\right)\right),\tag
2.14
$$
where the notation ``$\sim$" means that the expressions are equal
up to  terms having lower total degree, and where for $i,j,k,\ell$
fixed, $A_{k\ell}^{ij}$ and $B_{k\ell}^{ij}$ are $N\times N$
matrix valued functions.  The term $u^{(1)}$ is obtained by
considering the terms in (2.4)-(2.9) of the form
$(S^2-1/4)\xi_j\xi_k\xi_p\xi_q$. We define the {\it coefficient
symbols} of $\Delta'$ by replacing each occurrence of $\partial_j$
in (2.14) by $\xi_j$.
$$
\s^{(2)}\ =\ \sum_{i,j}h_{ij}\xi_i\xi_j, \qquad\qquad \s^{(1)}\ =\
\sum_{i,j,k,\ell}A^{ij}_{k\ell}\ \xi_k\xi_\ell\ h_{ij}\qquad\qquad
\s^{(0)}\ =\ \sum_{i,j,k,\ell}B^{ij}_{k\ell}\ \xi_k\xi_\ell\
h_{ij}.\tag 2.15
$$
These are $N\times N$ matrix valued.  We also define the {\it
coefficient components} of $\Delta'$ by
$$
\s^{(1)}_{k\ell}\ =\ \sum_{i,j}A^{ij}_{k\ell} h_{ij}\qquad\qquad
\s^{(0)}_{k\ell}\ =\ \sum_{i,j}B^{ij}_{k\ell} h_{ij}.\tag 2.16
$$
 Then   $u_s^{(1)}$ is obtained by collecting all multiples of
 $(S^2-1/4)\xi_j\xi_k\xi_p\xi_q$ from (2.4)--(2.9), and its
 quadratic form is given by
$$
\langle h\,,\,(u_s^{(1)}(x,\xi)h\rangle \ :=\
(S^2-1/4)|\xi|^{n-2s-4} \trace\left(\left(
\s^{(2)}-2\s^{(1)}+4\s^{(0)} \right)^2\right).\tag 2.17
$$
The symbol $u_s^{(2)}$ is obtained from  (2.6) after subtracting
the term which goes into $u_s^{(1)}$. Its quadratic form is given
by
$$
\langle h\,,\,(u_s^{(2)}(x,\xi)h\rangle  = (2S-1) |\xi|^{n-2s-4}
\trace\left(\left( \s^{(1)} \right)^2\ -\ |\xi|^2\sum_{j}\left(
\sum_i \xi_i\ \s^{(1)}_{ij} \right)^2\right).\tag 2.18
$$
The term $u_s^{(3)}$ is computed from (2.7) and (2.8) and
represents the interaction of $\sum_{i,j}h_{ij}\partial^2_{ij}$
with the other terms after subtracting the terms which go into
$u_s^{(1)}$:
$$
\align \langle h\,,\,(u_s^{(3)}(x,\xi)h\rangle \ =\
(2S-1)|\xi|^{n-2s-4}&\Biggl((\xhx)
\trace \left(-\s^{(1)} -2\s^{(0)} \right) \\
&+\  |\xi|^2(\tr h)
\trace\left(-\s^{(1)} +2\s^{(0)} \right)\\
&+\ |\xi|^2 \trace \sum_{j}  \s^{(1)}_{ij}\
\left(\xi_i\,(\hx)_j+\xi_j\,(\hx)_i\right)\Biggr).\tag 2.19
\endalign
$$
The term $u_s^{(4)}$ is computed from (2.9) and represents the
interaction of $\sum_{i,j}h_{ij}\partial^2_{ij}$ with itself after
$u_s^{(1)}$ has been taken into account:
$$
\multline \langle h\,,\, (u_s^{(4)}(x,\xi)h\rangle \ =\
\dim E\ |\xi|^{n-2s-4}\\
 \times\ \left( \frac12 |\xi|^4|h|^2\ -\ 2S|\xi|^2 |\hx|^2\ +\
 (S+1/4)(\xhx)^2
 \ +\ (S-1)|\xi|^2(\tr h)(\xhx)\ +\ \frac14 |\xi|^4(\tr h)^2
  \right)\\
  =\ \dim E\ |\xi|^{n-2s-4}(S-1/2)\left( -2|\xi|^2|\hx|^2+(\xhx)^2
  +|\xi|^2(\xhx)(\tr h)\right)
  \\ +\ \dim E\ |\xi|^{n-2s}\left( \frac12\tr\left((H\Pi_\xi^\perp)^2\right)
  \ +\ \frac14\left(\tr H\Pi_\xi^\perp\right)^2\right).
\endmultline
\tag 2.20
$$
The fact that (2.11) holds may be verified by explicitly computing
$\langle h\,,\,u_s(x,\xi)h\rangle$ by substituting  (2.14) into
(2.3).
\medskip\medskip

\noindent{\bf 3. Calculation of $\Delta_p'$}.
\medskip

To apply these formulas to  the Bochner and de Rham Laplacians
$\Delta_p$, we must compute the terms of  $\Delta_p '$ of total
degree $2$ in the form (2.14). Let $(x^1,\dots,x^n)$ be local
coordinates on $M$, and let $\Lambda^p(M)$ denote the bundle of
alternating $(0,p)$ tensors on $M$.  Then the subsets
$I\subset\{1,2,\dots,n\}$ of size $p$ index a basis for
$\Lambda^p_x(M)$. The $I$th basis element is $dx^I=
dx^{i_1}\wedge\dots\wedge dx^{i_p}$, where $i_1,\dots,i_p$ are
chosen so that
 $i_1<i_2<\dots<i_p$ and $I=\{i_1,\dots,i_p\}$.
 For a $p$-form $F $, we write in
components
$$
F \ =\ \sum_{i_1,\dots,i_p=1}^n  F_{i_1 \dots i_p}
dx^{i_1}\otimes\dots\otimes dx^{i_p}\ =\  \sum_{i_1<\dots<i_p} F
_{i_1 \dots i_p} dx^{i_1}\wedge\dots\wedge dx^{i_p} \ =\ \sum_I
F_Idx^I.
$$
where $F_I:=F_{i_1\dots i_p}$, with $i_1,\dots,i_p$ defined so
that $i_1<\dots<i_p$ and $I=\{i_1,\dots,i_p\}$.  It is convenient
to write $\delta$ for $d^*$.

 \proclaim{Lemma 3.1} In coordinates which are normal at a
point $x$,
$$
\align & ((\nabla^*\nabla)'F)_{i_1\dots i_p} \ \sim\ \sum_{j,k}
h_{jk}\
\partial^2_{jk} F_{i_1\dots i_p} \ +\ \sum_{j,k}
\left( \partial_j h_{jk} - \frac12\partial_k h_{jj} \right)
\partial_{k} F_{i_1 \dots i_p}  \\
&+\ \sum_{j,k}\sum_{s=1}^p (-1)^{s+1}
 \ \left(\partial_{{i_s}}
h_{jk}+\partial_j h_{i_s\,k}-\partial_k h_{i_s\,j}\right)
\partial_j F_{k \,i_1\dots\widehat{i_s} \dots i_p} \\
&\ + \  \frac12\sum_{j,k}\sum_{s=1}^p
(-1)^{s+1}\left(\partial^2_{j\,i_s}h_{jk}-\partial^2_{jk}h_{i_s
j}+
\partial^2_{jj} h_{i_s k}\right) F_{k\,i_1 \dots\widehat{i_s}\dots
i_p}, \tag 3.1
\endalign
$$
and
$$
\align ((\delta' d  +d\delta')F )_{i_1 \dots i_p}\ &\sim\ \
\sum_{j,k} h_{jk}\
\partial^2_{jk} F _{i_1\dots i_p} \ +\ \sum_{j,k}
\left( \partial_j h_{jk} - \frac12\partial_k h_{jj} \right)
\partial_{k} F _{i_1 \dots i_p}  \\
&+\ \sum_{j,k}\sum_{s=1}^p (-1)^{s+1}
 \ \left(\partial_{{i_s}}
h_{jk}+\partial_j h_{i_s\,k}-\partial_k h_{i_s\,j}\right)
\partial_j F _{k \,i_1\dots\widehat{i_s} \dots i_p} \\
&\ + \  \sum_{j,k}\sum_{s=1}^p
(-1)^{s+1}\left(\partial^2_{i_s\,j}h_{jk}-\frac12\partial^2_{i_s\,k}
h_{jj}\right) F _{k\,i_1 \dots\widehat{i_s}\dots i_p} \\
& +\ \sum_{j,k}\sum_{s=1}^p\sum_{t\neq s} (-1)^{s+t+1_{t>s}}
\left(\partial^2_{i_sj}\,h_{i_t k}\right) F_{j\,k\,i_1\dots
\widehat{i_t}\dots\widehat{i_s} \dots i_p}. \tag 3.2
\endalign
$$
\endproclaim
\demo{Proof} To carry out this calculation, it is convenient to
use the Einstein summation convention that a repeated index in the
upper and lower position is summed over. Then
$$
\nabla^*\nabla F_{i_1 \dots i_p}\ =\ -F_{i_1\dots
i_p;j}^{\hphantom{i_1\dots i_p;j}j}\ =\ -g^{j\ell}\partial_\ell
F_{i_1 \dots i_s;j}\ +\ g^{j\ell}\sum_{s=1}^p (-1)^{s+1}
\Gamma_{\ell i_s}^k F_{k i_1 \dots \widehat{i_s} \dots i_p;j}\ +\
g^{j\ell} \Gamma_{\ell j}^k F_{i_1 \dots i_s;k}\tag 3.3
$$
where the Christoffel symbols $\Gamma_{ij}^k$ are defined by
$$
\Gamma_{ij}^k\ =\ \frac{g^{k\ell}}2\left( \partial_i
g_{j\ell}+\partial_j g_{i\ell}-\partial_\ell g_{ij}\right). \tag
3.4
$$
Now using the fact that ${\Gamma_{ij}^k}'$ is tensorial in $h$,
and computing in normal coordinates for $g$, we see that
$$
G^{j k}_i\ :=\ g^{j\ell} ({\Gamma_{\ell i}^k}')\ =\ \frac12\left(
 h^{k\ j}_{i;}\ +\
h^{jk}_{\hphantom{jk};i}\ -\ h^{j\ k}_{i;}\right),\tag 3.5
$$
and so
$$
\multline ((\nabla^*\nabla)' F)_{i_1 \dots i_p}\\ \sim \
h^{j\ell}\partial_\ell F_{i_1 \dots i_s;j}\ -\ g^{j\ell}
\partial_\ell({F_{i_1 \dots i_s;j}}') \ +\ \sum_{s=1}^p (-1)^{s+1}
G_{i_s}^{j k} \partial_j F_{k i_1 \dots \widehat{i_s} \dots
i_p}\ +\ G_{j}^{j k} \partial_k F_{i_1 \dots i_p}\\
\sim \ h^{j\ell}\partial_\ell F_{i_1 \dots i_s;j}\ +\ g^{j\ell}
\partial_\ell\sum_{s=1}^p (-1)^{s+1} {\Gamma_{ji_s}^k}'{F_{k i_1 \dots\widehat{i_s}
\dots i_p}} \ +\ \sum_{s=1}^p (-1)^{s+1} G_{i_s}^{j k}
\partial_j F_{k i_1 \dots \widehat{i_s} \dots i_p}\ +\ G_{j}^{j
k} \partial_k F_{i_1 \dots i_s}
\\
\sim \ h^{j\ell}\partial_\ell\partial_j F_{i_1 \dots i_s}\ +\
2\sum_{s=1}^p (-1)^{s+1} G_{i_s}^{j k}
\partial_j F_{k i_1 \dots \widehat{i_s} \dots i_p}\ +\ G_{j}^{j
k} \partial_k F_{i_1\dots i_s}\ +\ \sum_{s=1}^p (-1)^{s+1}
(\partial_j G_{i_s}^{j k}){F_{k i_1\dots\widehat{i_s} \dots i_p}}
\endmultline
$$
Substituting the formula for $G$ in (3.5) gives (3.1).

For the de Rham Laplacian, we have
$$
(dF)_{i_1 \dots i_{p+1}}\ =\  \sum_{s=1}^{p+1} (-1)^{s+1}
\partial_{{i_s}} F_{i_1 \dots \widehat{i_s} \dots i_{p+1}},\tag
3.6
$$
and
$$
(\delta F)_{i_1 \dots i_{p-1}}\ =\ -  F_{j\, i_1 \dots
i_{p-1};}^{\hphantom{j \,i_1 \dots i_{p-1};}j}.\tag 3.7
$$

But
$$
\multline - F_{j\, i_1 \dots i_{p-1};}^{\hphantom{j \,i_1 \dots
i_{p-1};}j}\ =\ -g^{j\ell} F_{j\,i_1 \dots i_{p-1};\ell}\\ =\
g^{j\ell}\left(-\partial_{\ell}F_{j\,i_1 \dots i_{p-1}}\ +\
\Gamma_{\ell j}^k F_{k\,i_1 \dots i_{p-1}}\ +\ \sum_{s=1}^{p-1}
(-1)^{s+1}\,\Gamma_{\ell \,i_s}^k F_{j\,k\,i_1\dots\widehat{i_s}
\dots i_{p-1}}\right),
\endmultline
\tag 3.8
$$
and so from (3.5),
$$
\left(-F_{j\, i_1 \dots i_{p-1};}^{\hphantom{j \,i_1 \dots
i_{p-1};}j}\right)'\ =\ h^{jk}F_{j \,i_1 \dots i_{p-1};k} \ +\
G_{j}^{jk} F_{k\,i_1 \dots i_{p-1}}\ +\ \sum_{s=1}^{p-1}
(-1)^{s+1} G_{i_s}^{j k} F_{j\,k\,i_1\dots \widehat{i_s} \dots
i_{p-1}}. \tag 3.9
$$
Set
$$
1_{t>s}\ =\ \cases 1 &\qquad t>s\\ 0& \qquad t\leq s\endcases.\tag
3.10
$$
Then
$$
\align (\delta' dF)_{i_1 \dots i_p}\ &=\ \ \left(-(dF)_{j\, i_1
\dots
i_p;}^{\hphantom{j \,i_1 \dots i_p;}j}\right)'\\
&\sim\  h^{jk}\partial_k(dF)_{j \,i_1 \dots i_{p}}\ +\ G_{j}^{jk}
(dF)_{k\,i_1 \dots i_{p}}\ +\ \sum_{s=1}^{p} (-1)^{s+1} G_{i_s}^{j
k} (dF)_{j\,k\,i_1\dots \widehat{i_s} \dots
i_{p}}\\
&=\ \   h^{jk}\ \partial^2_{jk} F_{i_1\dots i_p} \ +\ \sum_{s=1}^p
(-1)^s  \ h^{jk}\
\partial^2_{k \,{i_s}} F_{j\,i_1\dots
\widehat{i_s}\dots i_p}
 \\
& \qquad +\ G_{j}^{jk}\
\partial_{k} F_{i_1 \dots i_p}
 \ +\ G_j^{jk}\ \sum_{s=1}^p(-1)^s\
\partial_{i_s}
F_{k\,i_1 \dots \widehat{i_s}\dots i_p} \\
& \qquad +\ \sum_{s=1}^p(-1)^{s+1} \left(G_{i_s}^{j
k}-G_{i_s}^{kj}\right)\
\partial_{j} F_{k\,i_1 \dots \widehat{i_s} \dots i_p}\\
&\qquad +\ \sum_{s=1}^p\sum_{t\neq s} (-1)^{s+t+1_{t>s}}
G_{i_s}^{j k}\
\partial_{i_t} F_{j\,k\,i_1  \dots \widehat{i_t}\dots \widehat{i_s}\dots
i_p}.\tag 3.11\\
\endalign
$$
On the other hand,
$$
\align (d\delta' & F)_{i_1\dots i_p}\ =\ \sum_{s=1}^p (-1)^{s+1} \
\partial_{{i_s}} (\delta'F)_{i_1 \dots \widehat{i_s} \dots
i_{p}}\\
&=\  \sum_{s=1}^p (-1)^{s+1} \
\partial_{{i_s}} \left(h^{jk}F_{j \,i_1\dots\widehat{i_s} \dots i_p;k}
 \ +\ G_{j}^{jk} F_{k\,i_1 \dots\widehat{i_s}\dots
i_p}\ +\ \sum_{t\neq s}(-1)^{t+1_{s>t}} G_{i_t}^{j k}
F_{j\,k\,i_1\dots \widehat{i_t}\dots\widehat{i_s} \dots
i_p}\right)\\ \allowdisplaybreak &\sim\  \sum_{s=1}^p (-1)^{s+1} \
\left(\partial_{{i_s}} h^{jk}\right)\partial_k F_{j
\,i_1\dots\widehat{i_s} \dots i_p} \ +\ \sum_{s=1}^p (-1)^{s+1} \
 h^{jk}\partial^2_{i_s k}F_{j \,i_1\dots\widehat{i_s} \dots i_p} \\
&\ +\  \sum_{s=1}^p
(-1)^{s+1}\left(\partial_{i_s}G_{j}^{jk}\right) F_{k\,i_1
\dots\widehat{i_s}\dots i_p} \ +\ \sum_{s=1}^p (-1)^{s+1}
G_{j}^{jk}\partial_{i_s} F_{k\,i_1 \dots\widehat{i_s}\dots
i_p}  \\
&\ +\ \sum_{s=1}^p\sum_{t\neq s} (-1)^{s+t+1_{t>s}}
\left(\partial_{i_s}G_{i_t}^{j k}\right) F_{j\,k\,i_1\dots
\widehat{i_t}\dots\widehat{i_s} \dots i_p}
 \ +\ \sum_{s=1}^p\sum_{t\neq s} (-1)^{s+t+1_{t>s}}G_{i_t}^{j k}
\partial_{i_s} F_{j\,k\,i_1\dots \widehat{i_t}\dots\widehat{i_s} \dots i_p}
\tag 3.12
\endalign
$$
We notice that the second terms on the last three lines of (3.12)
cancel with terms in (3.11), and  hence
$$
\align &((\delta' d  +d\delta')F)_{i_1 \dots i_p}\ \sim\ h^{jk}\
\partial^2_{jk} F_{i_1\dots i_p} \ +\ G_{j}^{jk}\ \partial_{k} F_{i_1 \dots i_p}
 \\
& \ +\ \sum_{s=1}^p (-1)^{s+1} \ \left(\partial_{{i_s}}
h^{jk}\right)\partial_k F_{j \,i_1\dots\widehat{i_s} \dots i_p}
 \ +\ \sum_{s=1}^p(-1)^{s+1} \left(G_{i_s}^{j
k}-G_{i_s}^{kj}\right)\
\partial_{j} F_{k\,i_1 \dots \widehat{i_s} \dots i_p} \\
&\ + \  \sum_{s=1}^p
(-1)^{s+1}\left(\partial_{i_s}G_{j}^{jk}\right) F_{k\,i_1
\dots\widehat{i_s}\dots i_p} \ +\ \sum_{s=1}^p\sum_{t\neq s}
(-1)^{s+t+1_{t>s}} \left(\partial_{i_s}G_{i_t}^{j k}\right)
F_{j\,k\,i_1\dots \widehat{i_t}\dots\widehat{i_s} \dots i_p}
\endalign
$$
which gives (3.2) on substituting in the formula for $G$ in (3.5).
\enddemo
\medskip
\medskip
\noindent{\bf 4. Matrix entries of $\Delta_p'$.} \medskip

  As we
already discussed at the beginning of Section 3, if we fix local
coordinates on $M$ about the point $x$, then  $dx^I$ forms a basis
for $\Lambda^p_x$, where $I$ ranges over   subsets  of
$\{1,2,\dots,n\}$ of size $p$.  A $p$-form is then identified with
a vector valued function $(F_I)$, and the Bochner or de Rham
Laplacian at $x$ can be identified with a matrix of (scalar)
differential operators.  The purpose of this section is to
identify these matrix entries, or rather the terms of total degree
$2$ in these entries (see (2.13)). We will start by considering
the Bochner Laplacian $\nabla^*\nabla$, whose  terms of total
degree $2$ are given in (3.1), which we rewrite here:
$$
\align & ((\nabla^*\nabla)'F)_{i_1\dots i_p} \ \sim\ \sum_{j,k}
h_{jk}\
\partial^2_{jk} F_{i_1\dots i_p} \ +\ \sum_{j,k}
\left( \partial_j h_{jk} - \frac12\partial_k h_{jj} \right)
\partial_{k} F_{i_1 \dots i_p}  \\
&+\ \sum_{j,k}\sum_{s=1}^p (-1)^{s+1}
 \ \left(\partial_{{i_s}}
h_{jk}+\partial_j h_{i_s\,k}-\partial_k h_{i_s\,j}\right)
\partial_j F_{k \,i_1\dots\widehat{i_s} \dots i_p} \\
&\ + \  \frac12\sum_{j,k}\sum_{s=1}^p
(-1)^{s+1}\left(\partial^2_{j\,i_s}h_{jk}-\partial^2_{jk}h_{i_s
j}+
\partial^2_{jj} h_{i_s k}\right) F_{k\,i_1 \dots\widehat{i_s}\dots
i_p}, \tag 4.1
\endalign
$$
Suppose that $I=\{i_1,\dots,i_p\}$ where $i_1<\dots<i_p$.  The
terms on the top line of (4.1) give $II$ matrix entries. The
summands appearing in the second and third lines give $KI$ matrix
entries, where $K=(I\setminus \{i_s\})\cup\{k\}$.  We get $K=I$
when $k=i_s$, so  altogether, the $II$ matrix entry of (4.1) is
$$
\sum_{j,k} h_{jk}\
\partial^2_{jk}  \ +\ \sum_{j,k}
\left( \partial_j h_{jk} - \frac12\partial_k h_{jj} \right)
\partial_{k} \ +\ \sum_{j}\sum_{i\in I}
 \ \left(\partial_j h_{ii}\right)
\partial_j \ + \  \frac12\sum_{j}\sum_{i\in I}
\left(
\partial^2_{jj} h_{ii}\right).\tag 4.2
$$
It is convenient to introduce the  notation
$$
\chi_I(j)\ =\ \cases 1 \quad & j\in I,\\ 0 \quad & j\notin
I\endcases, \tag 4.3
$$
$$
\sgn_I(j)\ =\ \cases 1 \qquad & j\in I,\\ -1 \quad & j\notin
I.\endcases \tag 4.4
$$
After simplifying (4.2), we get

\proclaim{Lemma 4.1} The term of total degree $2$  in the $II$
matrix entry of $\nabla^*\nabla$ is
$$
\sum_{i,j} \left(h_{ij}\
\partial^2_{ij}\  +\
\partial_i  h_{ij}\partial_j\right)\ +\
\frac12\sum_{i,j}\sgn_I(i)\partial_j h_{ii}\partial_j\ +\
\frac12\sum_{i,j}\chi_I(i) (\partial^2_{jj}
  h_{ii}).\tag 4.5
$$
\endproclaim

The only other non-vanishing matrix entries  of total degree $2$
in $\nabla^*\nabla$ are $KI$ terms where $K=(I\setminus
\{i_s\})\cup\{k\}$, and $k\notin I$. These can easily be read off
from (4.1), although one has to pay attention to the sign.
\definition{Definition 4.2} For $J\subset \{1,2,\dots,n\}$ and $i\in
\{1,2,\dots,n\}$, define
$$
N(i,J)\ =\ \#\{j\in J:j<i\}. \tag 4.6
$$
The point of this definition is  that
$$
dx^{\{i\}\cup J}\ =\ (-1)^{N(i,J)} \,dx^i\wedge dx^J. \tag 4.7
$$
\enddefinition

\proclaim{Lemma 4.3} Suppose that $J$ is a subset of
$\{1,\dots,n\}$ of size $p-1$, and $I=J\cup\{i\}$ and $K=J\cup
\{k\}$, where $i,k\notin J$ and $i\neq k$. Then the term of total
degree $2$ in the  $KI$ matrix entry  of $\nabla^*\nabla$ is
$$
 (-1)^{N(i,J)+N(k,J)}
 \ \sum_j \left(\left(\partial_{{i}}
h_{jk}+\partial_j h_{i\,k}-\partial_k h_{i\,j}\right)
\partial_j \ +\
\frac12\left(\partial^2_{ij}h_{jk}-\partial^2_{jk}h_{ij}
+\partial^2_{jj}h_{ik}\right)\right). \tag 4.8
$$
\endproclaim
\medskip

We can carry out the same analysis for the de Rham Laplacian in
(3.2).  The only new feature is a new type of term appearing in
the last line of (3.2), which  we rewrite here:
$$
\align ((\delta' d  +d\delta')F )_{i_1 \dots i_p}\ &\sim\ \
\dots\dots\dots\dots\dots \\
& +\ \sum_{s=1}^p\sum_{t\neq s}\sum_{k,\ell} (-1)^{s+t+1_{t>s}}
\left(\partial^2_{i_sk}\,h_{i_t \ell}\right) F_{k\,\ell\,i_1\dots
\widehat{i_t}\dots\widehat{i_s} \dots i_p}. \tag 4.9
\endalign
$$
\medskip

\noindent By writing $J=I\setminus\{i_s,i_t\}$ and
$L=J\cup\{k,\ell\}$ and noting that
$$
F_{{i_1}\dots{i_p}}\ =\ (-1)^{s+t+1_{t>s}}
F_{{i_s}\,{i_t}\,{i_1}\dots\widehat{i_s}\dots \widehat{i_t}\dots
{i_p}},
$$
we see that  $LI$ matrix entry in the summand of (4.9) is
$$
(-1)^{N(i_s,J)+N(i_t,J)+N(k,J)+N(\ell,J)} \sgn((i_t-i_s)(\ell-k))\
\left(\partial^2_{i_sk}\,h_{i_t \ell}\right).\tag 4.10
$$
 The case $L=I$
can occur in two ways: $k=i_s$ and $k=i_t$, or  $j=i_t$ and
$k=i_s$.  Hence after carrying out the  $k,\ell$ sum (but not the
$s,t$ sum),  the $II$ matrix entry of  (4.9)  is
$$
\left(\partial^2_{i_si_s}\,h_{i_t i_t}\right) \ -\
\left(\partial^2_{i_si_t}\,h_{i_t i_s}\right).\tag 4.11
$$
Hence the $II$ matrix entry of  (3.2) is
$$
\multline\sum_{j,k} h_{jk}\
\partial^2_{jk}  \ +\ \sum_{j,k}
\left( \partial_j h_{jk} - \frac12\partial_k h_{jj} \right)
\partial_{k}  \\
 + \  \sum_{i,j}\chi_I(i) \left((\partial_j h_{i\,i})\partial_j\
+\
\partial^2_{i\,j}h_{i j }-\frac12\partial^2_{i\,i}
h_{jj}\right)  \ +\ \sum_{i,j}\chi_I(i)\chi_{I'}(j)
\left(\partial^2_{ii}\,h_{jj}-\partial^2_{ij}h_{ij} \right).
\endmultline
\tag 4.12
$$
Simplifying this, we get the following.
 \proclaim{Lemma 4.4} The term of total degree $2$ in the $II$ matrix
entry of  $d^*d+dd^*$ is
$$
\align
 & \sum_{i,j} \left(h_{ij}\
\partial^2_{ij}\  +\
\partial_i  h_{ij}\partial_j\right)\ +\
\frac12\sum_{i,j}\sgn_I(j)\partial_i h_{jj}\partial_i\\
&\qquad\qquad\qquad\qquad\qquad\qquad \ +\
\sum_{i,j}\chi_I(i)\chi_{I'}(j) (\partial^2_{ij} h_{ij})\ +\
\frac12\sum_{i,j}\chi_I(i)\sgn_I(j)
\partial^2_{ii}
  h_{jj}.\tag 4.13
\endalign
$$
\endproclaim

 Now let $\bar J$ be a subset of $\{1,\dots,n\}$ of size
 $p-1$ and let $p,q$ with $p\neq q$ be such that
 $p\notin \bar J$ and $q\notin \bar J$.  Set $I=\bar J\cup\{p\}$ and
 $K=\bar J\cup\{q\}$. The  the summand in (4.9)
  gives  an $LI$ matrix entry, and
  $L=K$ if one of the following four cases holds:
 $$
 \align
 & i_s=p,\quad k=q,\quad i_t=\ell,\qquad\qquad\qquad\qquad
 i_s=p,\quad \ell=q,\quad i_t=k, \tag 4.14\\
 & i_t=p,\quad k=q,\quad i_s=\ell,\qquad\qquad\qquad\qquad
 i_t=p,\quad \ell=q,\quad i_s=k.
 \endalign
 $$
Hence the $KI$ matrix entry of (4.9) after summing is
$$
(-1)^{N(p,\bar J)+N(q,\bar J)} \sum_{j\in I}
\left(\partial^2_{pq}h_{jj}\ +\
\partial^2_{jj}\,h_{pq} \ -\
\partial_{jp}^2 h_{jq}\ -\ \partial_{jq}^2 h_{jp} \right) .\tag
4.15
$$
In order to see this, it is important to note that the sets $J$
from (4.10) and  $\bar J$ differ by one element, $i_s$ or $i_t$.
By considering the other terms in (3.2), we get the following.
\proclaim{Lemma 4.5}
 For $I=J\cup\{i\}$ and $K=J\cup
\{k\}$, where $i,k\notin J$ and $i\neq k$, the term of total
degree $2$ in the $KI$ matrix entry of $d^*d+dd^*$ is
$$
\align & (-1)^{N(i,J)+N(k,J)}\biggl(
 \ \sum_j \left(\left(\partial_{{i}}
h_{jk}+\partial_j h_{i\,k}-\partial_k h_{i\,j}\right)
\partial_j \ +\  (\partial^2_{i\,j}h_{jk})-\tfrac12(\partial^2_{i\,k}
h_{jj})\right)\\
&\qquad\qquad\qquad\qquad\qquad\qquad\qquad +\  \sum_{j}
\chi_I(j)\left(\partial^2_{ik}h_{jj}\ +\
\partial^2_{jj} h_{ik}\ -\  \partial^2_{ij}
h_{jk}\ -\ \partial^2_{jk}h_{ij}\right)\biggr) \tag 4.16\\
\endalign
$$
\endproclaim
Finally, similar considerations as in  the previous case give the
following. \proclaim{Lemma 4.6}
 Suppose $J$ is a subset of $\{1,\dots,n\}$ of size $p-2$,
 and  $I=J\cup\{i,j\}$, $L=J\cup
\{k,\ell\}$, where $i,j,k,\ell$ are distinct with $i<j$, $k<\ell$
and  $i,j,k,\ell\notin J$. Then the term of total degree $2$ in
the  $LI$ matrix entry  of $d^*d+dd^*$ is
$$
(-1)^{N(i,J)+N(j,J)+N(k,J)+N(\ell,J)}\ \left(\partial^2_{i
k}h_{j\ell}\ +\ \partial^2_{j\ell}h_{ik}\ -\
\partial^2_{i\ell}h_{jk}\ -\
\partial^2_{jk}h_{i\ell }\right).\tag 4.16
$$
Because $i,j,k,\ell$ are all distinct, we see that this vanishes
if $h$ is diagonal.
\endproclaim

\medskip
\medskip
\noindent{\bf 5. Proof of Theorem 1 for the Bochner Laplacian.}
\medskip

Combining the following Lemma with (2.12) and (2.20), gives
Theorem 1 for the Bochner Laplacian, see  (1.7), (1.8) and (1.9).
(Here in the proof we are dealing with the real case.  The factor
$2$ in (1.7) is the dimension of the complex line.)

 \proclaim{Lemma 5.1} For the Bochner Laplacian $\nabla^*\nabla$
on $p$ forms,
$$
\langle h\,,\,u_s^{(1)}(x,\xi)h \rangle\ =\ (S^2-1/4)|\xi|^{n-2s}\lm n\\
p\mr \left( \tr(H\Pi_\xi^\perp)\right)^2.\tag 5.1
$$
$$
\multline
\langle h\,,\,u^{(2)}(x,\xi)h\rangle \ =\ (S-1/2)|\xi|^{n-2s}4\lm n-2\\
p-1\mr \tr\left( (H\Pi_\xi^\perp)^2\right)\\ +\ (S-1/2)|\xi|^{n-2s-4} \lm n\\
p\mr \left( -2|\xi|^2|\hx|^2\ +\ 2(\xhx)^2 \right).
\endmultline
\tag 5.2
$$
$$
\multline \langle h\,,\, u^{(3)}(x,\xi)h\rangle \ =\
(S-1/2)|\xi|^{n-2s}
\lm n\\ p\mr \left( \tr (H\Pi_\xi^\perp)\right)^2\\ +\ (S-1/2)|\xi|^{n-2s-4} \lm n\\
p\mr \left( 4|\xi|^2|\hx|^2\ -\ 3(\xhx)^2 \ -\ |\xi|^2(\xhx)(\tr
h) \right).
\endmultline
\tag 5.3
$$
\endproclaim

\demo{Proof} For $j=0,1,2$, and subsets $I,K$ of $\{1,\dots,n\}$
of size $p$,  write $\s^{(j)KI}$ for the $KI$ matrix element of
$\s^{(j)}$, and write $\s^{(j)KI}_{pq}$ for the $KI$ matrix
element of $\s^{(j)}_{pq}$. Now  for a  general operator on
$p$-forms, we can expand the trace in (2.17) to get
$$
\multline \langle h\,,\,u_s^{(1)}h\rangle\ =\
(S^2-1/4)|\xi|^{n-2s-4} \Biggl(\sum_I\left(
\s^{(2)II}-2\s^{(1)II}+4\s^{(0)II} \right)^2 \\ +\ \sum_{K\neq
I}\left( \s^{(2)KI}-2\s^{(1)KI}+4\s^{(0)KI} \right)\left(
\s^{(2)IK}-2\s^{(1)IK}+4\s^{(0)IK} \right)\Biggr).
\endmultline
\tag 5.4
$$
For the Bochner Laplacian, we see from Lemma 4.1, we see that
$$
\s^{(2)II}\ =\ \xhx,\qquad\qquad \s^{(1)II}\ =\ (\xhx)\ +\
\frac12|\xi|^2\sum_i \sgn_I(i)h_{ii}, \qquad\qquad \s^{(0)II}\ =\
\frac12|\xi|^2\sum_i\chi_I(i)h_{ii}. \tag 5.5
$$
Hence
$$
\s^{(2)II}\ -\ 2\s^{(1)II}\ +\ 4\s^{(0)II} \ =\ -(\xhx)\
 +\  |\xi|^2(\tr h).\tag 5.6
$$
On the other hand, let $I,J,K$ be as in Lemma 4.3. Then
$$
\align &\s^{(1)KI}\ =\ (-1)^{N(i,J)+N(k,J)}\left(
\xi_i(\hx)_k-\xi_k(\hx)_i+|\xi|^2h_{ik}\right),\tag 5.7\\
& \s^{(0)KI}\ =\ (-1)^{N(i,J)+N(k,J)}\frac12\left( \xi_i(\hx)_k \
-\  \xi_k(\hx)_i\ +\ |\xi|^2 h_{ik}\right).
\endalign
$$
We see that
$$
\s^{(2)KI}\ -\ 2\s^{(1)KI}\ +\ 4\s^{(0)KI} \ =\ 0. \tag 5.8
$$
Hence from (5.6) and (5.8), we see that for the Bochner Laplacian,
$$
\multline \langle h\,,\, u_s^{(1)}h\rangle\ =\
(S^2-1/4)|\xi|^{n-2s-4}\sum_I \left(-(\xhx)\
 +\  |\xi|^2(\tr h)\right)^2\\ =\ (S^2-1/4)|\xi|^{n-2s-4}\lm n\\ p\mr \left(-(\xhx)\
 +\  |\xi|^2(\tr h)\right)^2,
 \endmultline \tag 5.9
$$
which gives (5.1).

Now we prove (5.2).  From (2.18), we can write
$$
\multline \langle h\,,\,u_s^{(2)}h\rangle  \ =\  (2S-1)
|\xi|^{n-2s-4} \Biggl(\sum_I\left(\left( \s^{(1)II} \right)^2\
-\ |\xi|^2\sum_{j}\left( \sum_i \xi_i\ \s^{(1)II}_{ij} \right)^2\right)\\
+\  \sum_{I\neq K}\left( \s^{(1)KI}\,\s^{(1)IK} \ -\
|\xi|^2\sum_{q}\left( \sum_p \xi_p\ \s^{(1)KI}_{pq} \right)
\left(\sum_r\xi_r\ \s^{(1)IK}_{rq} \right) \right)\Biggr).
\endmultline
\tag 5.10
$$
Now from (5.5) and (5.7) we have
$$
\left( \s^{(1)II} \right)^2\ =\ (\xhx)^2\ +\
|\xi|^2(\xhx)\sum_i\sgn_I(i) h_{ii}\ +\ \frac14|\xi|^4\left(
\sum_i\sgn_I(i) h_{ii}\right)^2.\tag 5.11
$$
and with $I,K$ as in Lemma 4.3,
$$
\multline \s^{(1)KI}\,\s^{(1)IK} \ =\ \left( \xi_i(\hx)_k \ -\
\xi_k(\hx)_i\ +\ |\xi|^2 h_{ik}\right)\left( \xi_k(\hx)_i \ -\
\xi_i(\hx)_k\ +\ |\xi|^2 h_{ik}\right)\\
=\ -\left( \xi_i(\hx)_k \ -\ \xi_k(\hx)_i\right)^2\ +\ |\xi|^4
h_{ik}^2.
\endmultline
\tag 5.12
$$

 With $I,J,K$ as in Lemma 4.3, we see from
Lemmas 4.1 and 4.3, that
$$
\align \s^{(1)II}_{pq}\ &=\
 h_{pq}\ +\ \frac{\delta_{pq}}2\sum_i\sgn_I(i) h_{ii},
 \qquad\qquad
  \sum \xi_p \s^{(1)II}_{pq}\ =\ (\hx)_q\ +\
\frac{1}2\,\xi_q\left(\sum_i\sgn_I(i) h_{ii}\right),\tag 5.13
\\
 \s^{(1)KI}_{pq}\ &=\  \delta_{ip} h_{qk}\ -\ \delta_{kp}h_{iq}\
+\ \delta_{pq} h_{ik},\qquad\qquad \sum \xi_p \s^{(1)KI}_{pq}\ =\
 \xi_i h_{qk}\ -\  \xi_k h_{iq}\ +\ \xi_q h_{ik}.
 \endalign
 $$
We see that
$$
|\xi|^2\sum_{j}\left( \sum_i \xi_i\ \s^{(1)II}_{ij} \right)^2 \ =\
|\xi|^2|\hx|^2\ +\ |\xi|^2(\xhx)\sum_i\sgn_I(i) h_{ii}\ +\
\frac14|\xi|^4\left( \sum_i\sgn_I(i) h_{ii}\right)^2, \tag 5.14
$$
and
$$
|\xi|^2\sum_{q}\left( \sum_p \xi_p\ \s^{(1)KI}_{pq} \right)
\left(\sum_r\xi_r\ \s^{(1)IK}_{rq} \right) \ =\ -|\xi|^2\sum_q
(\xi_ih_{kq}-\xi_k h_{iq})^2\ +\ |\xi|^4h_{ik}^2. \tag 5.15
$$
Putting (5.11), (5.12), (5.14) and (5.15) into (5.10), we get
$$
\multline \langle h\,,\, u^{(2)}h\rangle\ =\
 |\xi|^{n-2s-4} (2S-1)\ \Biggl(\sum_I
  \left((\xhx)^2-|\xi|^2|\hx|^2 \right)
  \\ -\ \sum_I\sum_{i,k}\chi_I(i)\chi_{I'}(k)\left(
 \left( \xi_i(\hx)_k \ -\ \xi_k(\hx)_i\right)^2
 \ -\ |\xi|^2\sum_q (\xi_ih_{kq}-\xi_k h_{iq})^2\right)
 \Biggr).
\endmultline
\tag 5.16
$$
Here, we use the fact that we can sum over all  $I$ and $K$ of the
form in Lemma 4.3, by summing over $I$, and over those $i,k$ with
$i\in I$ and $k\notin I$.  Using (7.1), we find that (5.16) equals
$$
\multline |\xi|^{n-2s-4} (2S-1)\ \Biggl(\lm n\\ p\mr
  \left((\xhx)^2-|\xi|^2|\hx|^2 \right)
  \\ +\ \lm n-2\\ p-1\mr \sum_{i,k}
  \left(|\xi|^2\sum_q (\xi_ih_{kq}-\xi_k h_{iq})^2
\ - \ \left( \xi_i(\hx)_k \ -\ \xi_k(\hx)_i\right)^2 \right)
 \Biggr).
\endmultline
\tag 5.17
$$
This equals (5.2).
\medskip
Now from (2.19),
$$
\align \langle h\,,\,u_s^{(3)}h\rangle \ &=\ (2S-1)|\xi|^{n-2s-4}
\sum_I \Biggl((\xhx)  \left(-\s^{(1)II} -2\s^{(0)II} \right) \\
&\qquad+\  |\xi|^2(\tr h) \left(-\s^{(1)II} +2\s^{(0)II} \right)\
+\ |\xi|^2  \sum_{j} \s^{(1)II}_{ij}\
\left(\xi_i\,(\hx)_j+\xi_j\,(\hx)_i\right)\Biggr) \tag 5.18
\endalign
$$
Using (5.5) and (5.13),  for the Bochner Laplacian we get
$$
\align & -\s^{(1)II} -2\s^{(0)II}\ =\  -(\xhx)\ -\ |\xi|^2
\sum_{i} \left( \frac12+ \sgn_I(i)\right) h_{ii} \\
&-\s^{(1)II} +2\s^{(0)II}\ =\  -(\xhx)\ +\ \frac12|\xi|^2(\tr h)  \\
& \sum_{j} \s^{(1)II}_{ij}\
\left(\xi_i\,(\hx)_j+\xi_j\,(\hx)_i\right)\ =\  2|\hx|^2\ +\
(\xhx)\sum_i \sgn_I(i) h_{ii}.
\endalign
$$
Hence putting this into (5.18),
$$
\align
 \langle h\,,\,u_s^{(3)}h\rangle \ & =\
(2S-1)|\xi|^{n-2s-4} \sum_I \Biggl((\xhx)  \left( -(\xhx)\ -\
|\xi|^2
\sum_{i} \left( \frac12+ \sgn_I(i)\right)h_{ii}   \right) \\
&\qquad+\  |\xi|^2(\tr h) \left(  -(\xhx)\ +\ \frac12|\xi|^2(\tr
h)\right)\ +\ |\xi|^2  \left(2|\hx|^2\ +\ (\xhx)\sum_i \sgn_I(i)
h_{ii} \right)\Biggr)
\endalign
$$
The terms involving $\sum_i \sgn_I(i) h_{ii}$ cancel, and we are
left with
$$
(2S-1)|\xi|^{n-2s-4}\lm n\\ p\mr\left( 2|\xi|^2|\hx|^2\ -\
(\xhx)^2\ -\ \frac32 |\xi|^2( \xhx)(\tr h)\ +\ \frac12|\xi|^4(\tr
h)^2\right),
$$
which equals (5.3).
\enddemo
\medskip
\medskip

\noindent{\bf 6. Proof of Theorem 1 for the de Rham Laplacian.}
\medskip

Combining the following Lemma with (2.12) and (2.20), gives
Theorem 1 for the de Rham Laplacian, see (1.7), (1.10) and (1.11).
(Here in the proof we are dealing with the real case.  The factor
$2$ in (1.7) is the dimension of the complex line.)

\proclaim{Lemma 5.2} For the de Rham Laplacian on $p$ forms,
$$
\multline
\langle h\,,\, u^{(1)}_sh\rangle \\ =\ (S^2-1/4)|\xi|^{n-2s}\left(4\lm n-2\\
p-1\mr \tr\left( (H\Pi_\xi^\perp)^2\right)\ +\ \left( \lm n\\ p\mr
-4\lm n-2\\ p-1\mr\right) \left( \tr
(H\Pi_\xi^\perp)\right)^2\right).
\endmultline
$$
$$
\multline
\langle h\,,\, u^{(2)}_sh\rangle \ =\ (S-1/2)|\xi|^{n-2s}4\lm n-2\\
p-1\mr \tr\left( (H\Pi_\xi^\perp)^2\right)\\ +\ (S-1/2)|\xi|^{n-2s-4} \lm n\\
p\mr \left( -2|\xi|^2|\hx|^2\ +\ 2(\xhx)^2 \right).
\endmultline
$$
$$
\multline \langle h\,,\, u^{(3)}_sh\rangle \ =\
(S-1/2)|\xi|^{n-2s}
\left( \lm n\\ p\mr-4\lm n-2\\
p-1\mr\right) \left( \tr (H\Pi_\xi^\perp)\right)^2\\ +\ (S-1/2)|\xi|^{n-2s-4} \lm n\\
p\mr \left( 4|\xi|^2|\hx|^2\ -\ 3(\xhx)^2 \ -\ |\xi|^2(\xhx)(\tr
h) \right).
\endmultline
$$
\endproclaim

\demo{Proof} We   work in coordinates $\xi$ in which $h$ is
diagonal. This is essential for  the calculation of $u^{(1)}_s$,
which becomes extremely lengthy if we drop this assumption. We
know that $u^{(1)}_s$ is itself independent of the coordinates,
because it is the coefficient of  $u_s$ in the variable $s^2$.  To
compute $u^{(1)}_s$ we  apply (5.4) and use the notation given
there. By Lemma 4.4, for the de Rham Laplacian,
$$
\align & \s^{(2)II}\ =\ \xhx,\\
& \s^{(1)II}\ =\ (\xhx)\
+\ \frac12|\xi|^2\sum_i \sgn_I(i)h_{ii}, \\
& \s^{(0)II}\ =\ \sum_{i,j}\chi_I(i)\chi_{I'}(j) \xi_i\xi_j
h_{ij}\ +\ \frac12\sum_{i,j}\chi_I(j)\sgn_I(i) \xi_j^2  h_{ii} \
=\ \frac12\sum_{i,j}\chi_I(j)\sgn_I(i) \xi_j^2  h_{ii}, \tag 6.4
\endalign
$$
where the last equality holds since $h$ is diagonal.
 Hence
$$
\s^{(2)II}\ -\ 2\s^{(1)II}\ +\ 4\s^{(0)II} \ =\ -(\xhx)  \  +\
\sum_{i,j} \sgn_I(i)\sgn_I(j) \xi_i^2 h_{jj}.\tag 6.5
$$
On the other hand, let $I,J,K$ be as in Lemma 4.5. Then
$\s^{(2)KI}=0$, and
$$
\align &\s^{(1)KI}\ =\ (-1)^{N(i,J)+N(k,J)}\left(
\xi_i(\hx)_k-\xi_k(\hx)_i+|\xi|^2h_{ik}\right),\tag 5.7\\
& \s^{(0)KI}\ =\ (-1)^{N(i,J)+N(k,J)}\Biggl( \xi_i(\hx)_k \ -\
\frac12 \xi_i\xi_k(\tr h)\\
&\qquad\qquad\qquad\qquad\qquad\qquad\qquad\qquad +\     \sum_{j}
\chi_I(j)\left(\xi_i\xi_k h_{jj}\ +\ \xi_j^2 h_{ik}\ -\ \xi_i\xi_j
h_{jk}\ -\ \xi_i\xi_kh_{ij}\right)\biggr)\Biggr). \tag 6.6
\endalign
$$
We see that
$$
\multline \s^{(2)KI}\ -\ 2\s^{(1)KI}\ +\ 4\s^{(0)KI} \\ =\
(-1)^{N(i,J)+N(k,J)} 2 \sum_{j} \sgn_I(j)\left(\xi_i\xi_k h_{jj}\
+\ \xi_j^2 h_{ik}\ -\ \xi_i\xi_j h_{jk}\ -\
\xi_i\xi_kh_{ij}\right)\\
 =\
(-1)^{N(i,J)+N(k,J)} 2\left(\left(\sum_j\sgn_I(j)h_{jj}\right)\ +\
h_{kk}\ -\ h_{ii}\right)\xi_i\xi_k,
\endmultline \tag 6.7
$$
where the last line follows because  $h$ is diagonal and $i\in I$,
$k\notin I$.  Because of Lemma 4.6, we have now considered all the
non-zero matrix entries of the coefficient symbols.  From (5.4) we
get
$$
\align &\langle h\,,\, u^{(1)}_sh\rangle \ =\
(S^2-1/4)|\xi|^{n-2s-4}\sum_I \\ &\left(\! \left(
 \!-(\xhx) +\sum_{i,j}\sgn_I(i)\sgn_I(j) \xi_i^2 h_{jj}\right)^2
\!\!  +  4\sum_{i,k} \chi_I(i)\chi_{I'}(k)\xi_i^2\xi_k^2
\left(\biggl(\sum_j\sgn_I(j)
 h_{jj}\biggr)^2 \!\! -  \left(h_{kk}-h_{ii}\right)^2\!\right)\!\right)
 \endalign
 $$
Multiplying out the brackets, we get
$$
\align \langle h\,,\, u^{(1)}_sh\rangle \ =\ (S^2- &
1/4)|\xi|^{n-2s-4}\Biggl(\lm n\\ p\mr (\xhx)^2
\ -\ 2(\xhx)\sum_{i,j} \sgn_I(i)\sgn_I(j)\xi_i^2 h_{jj}\\
&\ +\ \sum_{i,j,k,\ell} \sum_I\sgn_I(j)\sgn_I(\ell)\ \left(
\sgn_I(i)\sgn_I(k)+4\chi_I(i)\chi_{I'}(k)\right)\
\xi_i^2\xi_k^2 h_{jj} h_{\ell\ell}\\
&\ +\ 4\sum_{i,k}\sum_I\chi_I(i)\chi_{I'}(k)
\xi_i^2\xi_k^2\left(2h_{ii}h_{kk}- h_{ii}^2- h_{kk}^2\right)
\Biggr).\tag 6.8
\endalign
$$
But we can considerably simplify the second line of (6.8) because
$$
\multline \sum_{i,k}\left(
\sgn_I(i)\sgn_I(k)+4\chi_I(i)\chi_{I'}(k)\right)\ \xi_i^2\xi_k^2
\\
=\ \sum_{i,k}\left(
\sgn_I(i)\sgn_I(k)+2\chi_I(i)\chi_{I'}(k)+2\chi_{I'}(i)\chi_I(k)\right)\
\xi_i^2\xi_k^2 \ =\ \sum_{i,k}  \xi_i^2\xi_k^2\ =\ |\xi|^4.
\endmultline
$$
Hence we get
$$
\multline \langle h\,,\, u^{(1)}_sh\rangle \ =\ C(s)(S^2-
1/4)|\xi|^{n-2s-4}\Biggl(\lm n\\ p\mr (\xhx)^2
\ -\ 2(\xhx)\sum_{i,j} \sum_I\sgn_I(i)\sgn_I(j)\xi_i^2 h_{jj}\\
\ +\ |\xi|^4\sum_{j,\ell} \sum_I\sgn_I(j)\sgn_I(\ell) h_{jj}
h_{\ell\ell}\ +\ 4\sum_{i,k}\sum_I\chi_I(i)\chi_{I'}(k)
\xi_i^2\xi_k^2\left(2h_{ii}h_{kk}- h_{ii}^2- h_{kk}^2\right)
\Biggr).\endmultline\tag 6.9
$$
From (7.1) and (7.2),
$$
\align & \sum_{i,j} \sum_I\sgn_I(i)\sgn_I(j)\xi_i^2 h_{jj}
 \ =\ \left( \lm n\\ p\mr - 4\lm n-2\\ p-1\mr \right)|\xi|^2 (\tr h)\ +\
 4\lm n-2\\ p-1\mr (\xhx),\\
&\sum_I\sgn_I(j)\sgn_I(\ell) h_{jj} h_{\ell\ell}\ =\ \left( \lm
n\\ p\mr - 4\lm n-2\\ p-1\mr \right) (\tr h)^2\ +\
 4\lm n-2\\ p-1\mr |h|^2,\\
 &\sum_{i,k}\sum_I\chi_I(i)\chi_{I'}(k)
\xi_i^2\xi_k^2\left(2h_{ii}h_{kk}- h_{ii}^2- h_{kk}^2\right) \ =\
2\lm n-2\\ p-1\mr\left( (\xhx)^2\ -\ |\xi|^2|\hx|^2\right).
 \endalign
 $$
Hence putting this into (6.9), we get
$$
\multline \langle h\,,\, u^{(1)}_sh\rangle \ =\
(S^2-1/4)|\xi|^{n-2s-4} \biggl( \lm n\\
p\mr (\xhx)^2\ -\ 2\left( \lm n\\ p\mr-4\lm n-2\\
p-1\mr\right)|\xi|^2(\xhx)(\tr h)\\ +\  \left( \lm n\\ p\mr-4\lm n-2\\
p-1\mr\right)|\xi|^4 (\tr h)^2\ +\ 4\lm n-2\\
p-1\mr |\xi|^4|h|^2\ -\ 8\lm n-2\\ p-1\mr |\xi|^2|\hx|^2,
\endmultline
$$
which equals (6.1).

The proof of (6.2) is the same as the proof of (5.2), because
$\s^{(1)}$ is the same for the Bochner and Hodge Laplacians.

 Finally we compute $u_s^{(3)}$ using (5.18).  From  (6.4) and Lemma 4.4, we have
 $$
\align & -\s^{(1)II} -2\s^{(0)II}\ =\  -(\xhx)\ -\ |\xi|^2\sum_i
\sgn(i)h_{ii}\ -\ \frac12
\sum_{i,j}\sgn_I(i)\sgn_I(j) \xi_j^2 h_{ii} \\
&-\s^{(1)II} +2\s^{(0)II}\ =\  -(\xhx)\ +\ \frac12\sum_{i,j}
\sgn_I(i)\sgn_I(j) \xi_i^2
h_{jj}  \\
& \sum_{j} \s^{(1)II}_{ij}\
\left(\xi_i\,(\hx)_j+\xi_j\,(\hx)_i\right)\ =\  2|\hx|^2\ +\
(\xhx)\sum_i \sgn_I(i) h_{ii}. \tag 6.10
\endalign
$$
 and so putting this into (5.18), we get
$$
\multline \langle h\,,\,u_s^{(3)}h\rangle \ =\
(2S-1)|\xi|^{n-2s-4}\sum_I \\ \Biggl((\xhx) \left( -(\xhx)\ -\
|\xi|^2\sum_i \sgn_I(i)h_{ii}\ -\ \frac12
\sum_{i,j}\sgn_I(i)\sgn_I(j) \xi_j^2 h_{ii}
\right)\\
+\  |\xi|^2(\tr h) \left( -(\xhx)\ +\
\frac12\sum_{i,j}\sgn_I(i)\sgn_I(j) \xi_j^2 h_{ii}\right) \ +\
|\xi|^2  \left(2|\hx|^2\ +\ (\xhx)\sum_i \sgn_I(i) h_{ii}
\right)\Biggr).
\endmultline
$$
The two terms involving a single sign functions cancel, and we get
$$
\multline
 \langle h\,,\,u_s^{(3)}h\rangle \ =\ (2S-1)|\xi|^{n-2s-4}\Biggl(
\lm n\\ p\mr \left(-(\xhx)^2 \ -\ |\xi|^2(\xhx)(\tr h)\ +\
2|\xi|^2|\hx|^2\right)\\
+\ \frac12(|\xi|^2(\tr h)-(\xhx))\
\sum_{i,j}\sum_I\sgn_I(i)\sgn_I(j) \xi_j^2 h_{ii} \Biggr) .
\endmultline
$$
From (7.2), this gives
$$
\multline
 \langle h\,,\,u_s^{(3)}h\rangle \ =\ (2S-1)|\xi|^{n-2s-4}\Biggl(
\lm n\\ p\mr \left(-(\xhx)^2 \ -\ |\xi|^2(\xhx)(\tr h)\ +\
2|\xi|^2|\hx|^2\right)\\
+\ \frac12\left( \lm n\\ p\mr - 4\lm n-2\\ p-1\mr
\right)(|\xi|^2(\tr h)-(\xhx))|\xi|^2(\tr h)\ +\ 2\lm n-2\\ p-1\mr
(|\xi|^2(\tr h)-(\xhx))(\xhx).
\endmultline
$$
Collecting terms, we get
$$
\multline \langle h\,,\,u_s^{(3)}h\rangle \ =\
(2S-1)|\xi|^{n-2s-4}\Biggl(2\lm n\\ p\mr|\xi|^2 |\hx|^2
 \ +\  \left( -\lm n\\ p\mr-2 \lm n-2\\ p-1\mr
\right) \left(\xhx \right)^2\\
 +\ \left(-\frac32\lm n\\ p\mr +4\lm n-2\\ p-1\mr\right) |\xi|^2(\tr h)(\xhx)
 \ +\ \frac12\left( \lm n\\ p\mr-4\lm n-2\\ p-1\mr \right) |\xi|^4\left(
\tr h\right)^2\Biggr),
\endmultline
$$
which equals (6.3).
\medskip
\enddemo

\medskip
\medskip
\noindent{\bf 7. Free terms.}
\medskip
In Sections 5 and 6, we use the following elementary formulas. The
sums here are over all subsets $I$ of $\{1,\dots,n\}$ of size $p$.
$$
\sum_I \chi_I(i)\chi_{I'}(j)\ =\ \lm n-2 \\ p-1\mr (1-
\delta_{ij}).\tag 7.1
$$
$$
 \sum_{I} \sgn_I(j)\sgn_I(k)\ =\
 \left( \lm n\\ p\mr - 4\lm n-2\\ p-1\mr \right) \ +\ 4\lm n-2\\
p-1\mr \delta_{jk}.\tag 7.2
$$
Proving such formulas is simple. For example, for (7.1), one
considers separately the cases when $i=j$ and $i\neq j$.  When
$i\neq j$, the left hand side is the number of sets $I$ which
contain $i$ but not $j$, which is $\lm n-2\\ p-1\mr$. When $i=j$
the left hand side is clearly zero. For (7.2), we first expand
each sign function: $\sgn_I=\chi_I-\chi_{I'}$, and then deal which
each product of characteristic functions separately. \medskip

We end with a remark on another approach to carrying out the
calculation of $u_s$.  For the sums in (7.1) and (7.2), we refer
to the coefficient on the right hand side which does not involve
any delta functions as the {\it free term} of the sum. This is the
value of the sum when all the fixed variables are distinct. For
example, the free term of the sum on the left in (7.1) is $\lm n-2\\
p-1\mr$, and the free term of
$$
\sum_I \sgn_I(i)\sgn_I(j)\sgn_I(k)\sgn_I(\ell) \tag 7.3
 $$
 is
$16\lm n-4\\ p-2\mr - 8\lm n-2\\ p-1\mr+\lm n\\ p\mr$. If one
tries to carry out the calculation of $\langle h\,,\,u_sh\rangle$
for the de Rham Laplacian without using coordinates in which $h$
is diagonal, one is faced with several sums similar to (7.3). It
turns out that one can replace such sums by their free terms, and
the resulting answer is  correct. The idea behind this is that the
non-free terms always give rise to quadratic expressions in $h$
which are not coordinate independent, whereas $u_s$ is coordinate
independent.  If one could prove that  the coordinate dependent
terms arising from the non-free terms cannot add up in the end to
produce a non-vanishing coordinate free term, then this would give
another reasonably short way to carry out the calculation of
$u_s$.

%%%%%%%%%%%%%%%%%%%%%%%%%%%%%%%%%%%%%%%%%%%%%%%%%%%%%%%%%%%%%%%%%%%%%%%%%%%%%%%%%%%%%%%%%%
%%%%%%%%%%%%%%%%%%%%%%%%%%%%%%%%%%%%%%%%%%%%%%%%%%%%%%%%%%%%%%%%%%%

%%%%%%%%%%%%%%%%%%%%%%%%%%%%%%%%%%%%%%%%%%%%%%%%%%%%%%%%%%%%%%%%%%%%%%%%%%%%%%%%%%%%%%%%%%%%%
%%%%%%%%%%%%%%%%%%%%%%%%%%%%%%%%%%%%%%%%%%%%%%%%%%%%%%%%%%%%%%%%%%%%%%%%%%%%%%%%%%%%%%%%%%%%%

\vskip20pt

\centerline{References}

\roster

\item"[Ba]" G. Baker and J. Dodziuk,  Stability of spectra of
Hodge-de Rham Laplacians, {\it Math. Z.} {\bf 224} (1997), no. 3,
327--345.

\item"[BGV]"  N. Berline, E. Getzler, and M. Vergne, {\it Heat kernels
and Dirac operators}. Grundlehren der Mathematischen
Wissenschaften, 298. Springer-Verlag, Berlin, 1992.

\item"{[Br]}"
T.  Branson,  Sharp inequalities, the functional determinant,
   and the complementary series. {\it Trans. Amer. Math. Soc.}
   {\bf 347} (1995),  3671--3742.

\item"{[BCY]}"
T. Branson, S. Y. A. Chang and P. Yang, Estimates  and extremals
for zeta function determinants on four-manifolds, {\it Commun.
Math. Phys.} {\bf 149}, (1992), 241-262.

\item"{[B\O]}"
T. Branson and B. {\O}rsted, { Conformal geometry and global
invariants}, {\it Differential Geometry and its Applications} {\bf
1}, 279-308, (1991).

\item"[BFKM]"
D. Burghelea, L. Friedlander, T. Kappeler and P. McDonald, On the
functional logdet and related flows on the space of closed
embedded curves on $S\sp
   2$. {\it J. Funct. Anal.} {\bf 120}, (1994), 440--466.

\item"{[CY]}"
S. Y. A. Chang and P. Yang, Extremal metrics of zeta function
determinants on 4-manifolds, {\it Annals of Math.} {\bf 142}
(1995) 171-212.

\item"[CQ]"
S. Y. A. Chang and J. Qing,  The zeta functional determinants on
   manifolds with boundary.
    {\it J. Funct. Anal.} {\bf 147}
   (1997),  327--399.

\item
"[Ch]" J. Cheeger, Analytic torsion and the heat equation, {\it
Ann. Math.}, {\bf 109}, 259-322 (1979).

\item"[CT]"
S. Chanillo, F. Treves, On the Lowest Eigenvalue of the Hodge
Laplacian. {\it J. Diff. Geom.} (1997), {\bf 45}, 273-287.

\item
"[CC]" B. Colbois, G. Courtois, A Note on the First Nonzero
Eigenvalue of the Laplacian Acting on $p$-forms, {\it Manuscripta
Math.}, {\bf 68}, 143-160 (1990).

\item
"[Do]" J. Dodziuk, Eigenvalues of the Laplacian on Forms, {\it
Proc. AMS}, {\bf 85}, 437-443 (1982).

\item"[Fe]" H. Fegan,  The spectrum of the Laplacian on forms over a
Lie group. {\it Pacific J. Math.} {\bf 90} (1980),  373--387.

\item"[Go]"
R. Gornet, A New Construction of Isospectral Riemannian
nilmanifolds with Examples, {\it Michigan Math J.} {\bf 43}
(1996), no. 1, 159-188.

\item
"[HZ]" A. Hassell and S. Zelditch,  Determinants of Laplacians in
   exterior domains. {\it Internat. Math. Res. Notices},
   (1999),  {\bf 18}, 971--1004.

\item"[Lo]" J. Lott,  Collapsing and the differential form Laplacian: the
case of a smooth limit space. {\it Duke Math. J.} {\bf 114}
(2002),  267--306.

\item"[Mo1]"
C. Morpurgo, Sharp trace inequalities for intertwining operators
on $S\sp n$ and $R\sp n$. {\it Internat. Math. Res. Notices}
(1999), {\bf  20}, 1101--1117.

\item"[M\"{u}]" W. M\"{u}ller, analytic torsion and $R$-torsion of Riemannian
manifolds, {\it Adv. in Math.} {\bf 28} (1978), 233-305.

\item"[Ok1]"
K. Okikiolu, Critical metrics for the determinant of the Laplacian
in odd dimensions. {\it Ann. of Math.}
 {\bf 153} (2001),  471--531.

\item"[Ok2]"
K. Okikiolu, Hessians  of spectral zeta functions, Duke Math.
Journal, {\it to appear.}

\item"[Ok3]"
K. Okikiolu, Critical metrics for the zeta functions of a family
of scalar Laplacians, {\it Preprint.}

\item"{[OPS1]}"
B. Osgood, R. Phillips and P. Sarnak, Extremals of determinants of
Laplacians, {\it  J.  Funct. Anal.} {\bf  80} (1988), 148-211.

\item"{[OPS2]}"
B. Osgood, R. Phillips, and P. Sarnak,  Compact isospectral sets
of
   surfaces. {\it J. Funct. Anal.} {\bf 80} (1988),  212--234.

\item"{[OPS3]}"
B. Osgood, R. Phillips, and P. Sarnak,  Moduli space, heights and
   isospectral sets of plane domains. {\it Ann. of Math.} {\bf 129}
    (1989),  293--362.

\item"{[Po]}"
A. Polyakov, Quantum geometry of bosonic strings, {\it Phys. Lett.
B} {\bf 103} (1981), 207--210.

\item
"[RS]" D. Ray and I. Singer, Analytic torsion, {\it Proc.
Symposium Pure Math Vol. XXIII} AMS, Providence, R.I., 1973,
167-181.

\item
"[{ Ri}]" {K. Richardson}, {Critical points of the determinant of
the Laplace operator}, {\it Jour.  Funct. Anal.}, {\bf 122},
(1994) 52-83.

\item"{[Ro]}"
S. Rosenberg, The variation of the de Rham zeta function, {\it
Trans. AMS.} {\bf  299}, no. 2, (1987), 535-557.

\item"[Ta]" J. Takahashi Upper bounds for the eigenvalues of the
Laplacian on forms on certain Riemannian manifolds. {\it J. Math.
Sci. Univ. Tokyo} {\bf 6} (1999), no. 1, 87--99.

\endroster
\medskip
\medskip
\medskip

\centerline{Kate Okikiolu} \centerline{University of California,
San Diego} \centerline{okikiolu\@math.ucsd.edu}
\medskip

\centerline{Caitlin Wang} \centerline{University of California,
San Diego} \centerline{cwang\@math.ucsd.edu}

\end